\documentclass[a4paper,oneside,10.5pt]{article}%
\usepackage{amsmath}
\usepackage{amsfonts}
\usepackage{amssymb}
\usepackage{graphicx}%
\providecommand{\U}[1]{\protect \rule{.1in}{.1in}}

\newtheorem{theorem}{Theorem}[section]

\newtheorem{assumption}[theorem]{Assumption}
\newtheorem{example}[theorem]{Example}

\newtheorem{lemma}[theorem]{Lemma}

\newtheorem{remark}[theorem]{Remark}

\numberwithin{equation}{section}

\begin{document}

\title{ The necessary and sufficient conditions for stochastic differential systems with
multi-time states cost functional}
\author{Shuzhen Yang\thanks{Security Institute for Financial Studies, Shandong University,
Jinan, Shandong 250100, PR China. (yangsz@sdu.edu.cn). And Center for Mathematical Economics, Bielefeld University, Germany. }
\thanks{This work was supported by the Fundamental Research Funds of Shandong University (2015GN015); Supported by the China Scholarship Council (CSC);Supported by the Promotive research fund for excellent young and middle-aged scientists of Shandong Province (BS2015SF004); Supported by the China Postdoctoral Science
Foundation funded project (2015M570584).} }
\date{}
\maketitle

\textbf{Abstract}: From economics point of view,  we investigate a new optimal control problem driven by a stochastic differential equation with a multi-time states cost functional. By constructing a series of first-order adjoint equations, we establish the stochastic maximum principle and sufficient optimality conditions for this new optimal control problem. A constraints problem also be studied. In the end, we develop a near optimal control problem for a general cost functional.


\addcontentsline{toc}{section}{\hspace*{1.8em}Abstract}

\section{Introduction}
Let us consider the following optimal production planning problem with the uncertainties demand and the multi-time states cost functional  in the productive cycle. Notice that, the demand of the society is always uncertainties which could be described by a stochastic differential equation as follows,
$$
y(t)=y(0)+\displaystyle\int_0^tb(s)ds+\displaystyle\int_0^t\sigma(s)dW(s),
$$
where $W(\cdot)$ is a standard Brownian motion under a probability space $(\Omega,\mathcal{F},\{\mathcal{F}_t\}_{0\leq t\leq T},P)$. In order to meet the demand, the factory will change the production rate $u(\cdot)$ with the demand $y(\cdot)$, i.e.,
$$
X(t)=X(0)+\displaystyle\int_0^t\big{[}u(s)-y(s)\big{]}ds.
$$
Consider the production capacity, the production rate $u(\cdot)$ may satisfy
$$
0\leq u(s)\leq K,\ \ 0\leq s\leq T,
$$
where $K$ is a constant.

In fact, the factory may have some limitation for inventory level $X(\cdot)$ in the productive cycle, i.e., for $0\leq t_1\leq t_2\leq \cdots\leq t_n$, there are constraints for $(X(t_1),X(t_2),\cdots,X(t_n))$,
\begin{equation}
\label{incons}
0\leq E[X(t_i)]\leq \alpha_i,\ \ i=1,2,\cdots,n.
\end{equation}
On the other hand, the factory need to pay the running cost for inventory and production, we denote it as $f(X(t),u(t))$ at time $t\in[0,T]$. Also, the factory need pay the disposal cost for inventory level $X(\cdot)$ at time $(t_1,t_2,\cdots,t_n)$, in general, we denote it as $\Psi(X(t_1),X(t_2),\cdots,X(t_n))$. Thus, the cost functional as follows,
\begin{equation}
\label{incos-1}
J(u(\cdot))=E\big{[}\displaystyle\int_0^Tf(X(t),u(t))dt+\Psi(X(t_1),X(t_2),\cdots,X(t_n))\big{]}.
\end{equation}
In this study, we will consider the following general state processes with the cost functional (\ref{incos-1}),
\begin{equation}
\label{ine-1}
X(s)=\int_{0}^{s}b(X(t),u(t))dt+\int_{0}^{s}\sigma(X(t),u(t))dW(t).
\end{equation}
Also, some constraints conditions which similar with equation (\ref{incons}) is considered.

In the case  where $\Psi(X(t_1),X(t_2),\cdots,X(t_n))=\Psi(X(t_n))$ with $t_n=T$, there are many works concerning this subject. We refer Bensoussan \cite{B81} and Bismut \cite{B78} for the local maximum principle with the convex control set, and Peng \cite{P90}  for the global  maximum principle with general control domain which may not convex, for more see \cite{H76}. Recalling that dynamic programming with related HJB equations and maximum principle are powerful approaches for solving optimal control problems (see \cite{P90},\cite{P92}, \cite{WY08}, \cite{Y99} and \cite{P97}). The HJB equations derived for stochastic delay systems (see \cite{CW12}, \cite{M84} and \cite{M96}).

In our previous paper Yang \cite{Y16b}, the stochastic maximum principle for the above stochastic differential systems (\ref{ine-1}) with a general cost functional  is developed. Further, in \cite{Y16b}, the terminal cost functional is $\Psi(X_{[0,T]})$, where $X_{[0,T]}=X(s)_{0\leq s\leq T}$. However, there are some strong assumptions about Fr\'{e}chet
derivatives in \cite{Y16b}, and the structure of which is too complicity, for more details see \cite{GY16,Y16a}.

In this study, we will remove the assumptions about Fr\'{e}chet derivatives in \cite{Y16b}.  In the following, we present the details of this study. We first derive the maximum principle for the optimal control problem (\ref{incos-1}), the main difficult is that the cost functional has the part $\Psi(X(t_1),X(t_2),\cdots,X(t_n))$. At moment, for the limitation of technique, we assume that the control domain is convex. Then, we construct a series of first-order  adjoint equations and establish the stochastic maximum principle via the duality technique. In the following, we investigate the sufficient conditions for the optimal control problem (\ref{incos-1}). Motivate with the beginning of this section, a constrains problem  is developed and some usefull results is given.  To the best of our knowledge, the sufficient conditions and constrains problem for optimal control problem (\ref{incos-1}) are first investigated in this study.

The paper is organized as follows: In Section 2, we present the stochastic optimal
control problem and show some examples to describe our main results. The proof of maximum principle theorem is given in Section 3. In Section 4, the sufficient conditions for optimality problem is investigated. In addition, we develop the constrains problem and obtain some usefull results in Section 5. In the end, we investigate a near optimal control problem for dealing with the model in \cite{Y16b} via the results in Section 3.

\section{The optimal control problem}

Let $W$ be a $d$-dimensional standard Brownian motion defined on a complete
filtered probability space $(\Omega,\mathcal{F},P;\{ \mathcal{F}(t)\}_{t\geq
0})$, where $\{ \mathcal{F}(t)\}_{t\geq0}$ is the $P$-augmentation of the
natural filtration generated by the Brownian motion $W$.

Let $T>0$ be given, consider the following controlled stochastic differential
equation,
\begin{equation}
d{X}(s)=b(X{(s)},u(s))ds+\sigma (X{(s)},u(s))dW(s) ,\quad s\in(0,T],\label{ODE_1}%
\end{equation}
with the initial condition $X(0)=x$, where
$u(\cdot)=\{u(s),s\in \lbrack0,T]\}$ is a control process taking value in a
convex set $U$ of $\mathbb{R}^m$ and $b,\sigma$ are given deterministic functions.

In this study, we consider the following multi-time states cost functional, which related with different objective at different time.
\begin{equation}
J(u(\cdot))=%
E\big{[}{\displaystyle \int \limits_{0}^{T}}
f(X{(t)},u(t))dt+\Psi(X(t_1),\cdots,X(t_n))\big{]},\label{cost-1}%
\end{equation}
with  $0=t_0\leq t_1\leq \cdots \leq t_n=T$ and
\[%
\begin{array}
[c]{l}%
b:\mathbb{R}^m\times U\to \mathbb{R}^m,\\
\sigma:\mathbb{R}^{m}\times U\to \mathbb{R}^{m\times d},\\
f:\mathbb{R}^m\times U\to \mathbb{R},\\
\Psi:\mathbb{R}^{m\times n}\to \mathbb{R},\\
\end{array}
\]
we set $\sigma=(\sigma^1,\sigma^2,\cdots,\sigma^d)$, and $\sigma^j\in \mathbb{R}^m$ for $j=1,2,\cdots, d$.

Let $b,\sigma,f$ uniformly continuous and satisfy the following
linear growth and Lispschitz conditions.

\begin{assumption}
\label{ass-b}Suppose there exists a constant $c>0$ such that%
\[%
\begin{array}
[c]{c}%
\left| b(x_{1},u)-b(x_{2},u)\right| +\left| \sigma(x_{1},u)-\sigma(x_{2},u)\right|
 \leq c\left|x_1-x_2 \right|,\\
\end{array}
\]
$\forall(x_{1},u),(x_{2},u)\in{\mathbb{R}^m}\times U$.
\end{assumption}

\begin{assumption}
\label{assb-b2}Suppose there exists a constant $c>0$ such that
\[
\left|b(x,u)\right|+\left|\sigma(x,u)\right| \leq c(1+\mid x \mid),\quad \forall(x,u)\in{\mathbb{R}^m}\times U.
\]
\end{assumption}

\begin{assumption}
\label{ass-fai}Let $b,\sigma,f,\Psi$ be first differentiable at $x$ and $u$, and their derivatives in $x$ are continuous in $(x,u)$.
\end{assumption}

Let $\mathcal{U}[0,T]=\{u(\cdot)\in L^2_{\mathcal{F}}(0,T;U)\}.$ Suppose
Assumptions \ref{ass-b} and \ref{assb-b2} hold, then there exists a unique
solution $X$ for equation (\ref{ODE_1}) (see \cite{LS78}).

Minimize (\ref{cost-1}) over $ \mathcal{U}[0,T].$ Any $\bar{u}(\cdot)\in \mathcal{U}[0,T]$
satisfying
\begin{equation}
J(\bar{u}(\cdot))= \underset{u(\cdot)\in\mathcal{U}[0,T]}{\inf}J(u(\cdot)) \label{cost-2}%
\end{equation}
is called an optimal control. The corresponding state trajectory $(\bar{u}(\cdot),\bar{X}(\cdot))$ are called an optimal state trajectory and optimal pair.

At very beginning, we will show some examples to describe the main results in the following sections. The first example verifies the maximum principle (necessary conditions) for the cost functional (\ref{cost-1}), the second one describes the optimal production planning problem  with multi-time state constraints.
\begin{example}
Let $T=1$, the controlled stochastic differential equation as follows:
\begin{equation*}
d{X}^u(s)=u(s)dW(s),\text{ \  \ }  0 \leq s \leq 1,%
\end{equation*}
with the initial condition $X(0)=1$,  where
$u(\cdot)=\{u(s),0\leq s\leq 1\}$ is a control process taking values in a
compact set $U=[0,1]$.
The cost functional is
\begin{equation}
J(u(\cdot))=\underset{u\in \mathcal{U}[0,1]}{\inf}E[-2(X^u(\frac{1}{2}))^2+(X^u(1))^2] \label{2ex-11}%
\end{equation}
and we can verify that
 \begin{equation}
   ( \bar{u}(t),\bar{X}(t))=
   \begin{cases}
    (1,1+W(t)) &\mbox{if $0\leq t \leq \frac{1}{2}$},\\
  (0,1+W(\frac{1}{2}))&\mbox{if $ \frac{1}{2}  < t \leq 1$, }
   \end{cases}
  \end{equation}
is an optimal pair of systems (\ref{2ex-11}).

Next, we introduce the following  first-order adjoint  equations. \\
\begin{equation*}%
\begin{array}
[c]{ll}%
d{p}(t)= q(t)dW(t),\quad \frac{1}{2}  < t < 1,\\
p(1)= -2\bar{X}(1)
\end{array}
\label{2prin-11}%
\end{equation*}
 and
 \begin{equation*}%
 \begin{array}
[c]{ll}%
d{p}(t)= q(t)dW(t),\quad  0\leq t < \frac{1}{2}, \\
p(\frac{1}{2})= 4\bar{X}(\frac{1}{2})+p(\frac{1}{2}^{+}).
\end{array}
\label{2prin-12}%
\end{equation*}
The solutions of first-order adjoint equations as follows:
 \begin{equation}
   ( p(t),q(t))=
   \begin{cases}
  (2+2W(t),2)&\mbox{ $0\leq t < \frac{1}{2}$,}\\
  (-2-2W(\frac{1}{2}),0) &\mbox{ $ \frac{1}{2} < t \leq 1$. }\\
   \end{cases}
  \end{equation}
 Thus,
 \begin{equation}
   H_u(\bar{X}(t),\bar{u}(t),p(t),q(t))(v-\bar{u}(t))=
   \begin{cases}
   2(v-1) &\mbox{ $0\leq t < \frac{1}{2}$,}\\
   0 &\mbox{ $ \frac{1}{2} < t \leq 1$, }\\
        \end{cases}
  \end{equation}
 with $v\in [0,1]$, thus, the optimal control pair $(\bar{u}(\cdot),\bar{X}(\cdot))$ satisfies the Theorem \ref{Maximumprinciple}.

\end{example}

\begin{example}
In this example, we consider the optimal production planning problem which is given in Section 1. Let $T=1$, the controlled stochastic differential equation as follows:
\begin{equation}
\label{3exp-1}
\begin{array}
[c]{ll}
{X}^u(s)=\displaystyle\int_0^s\big{[}u(t)-y(t)\big{]}dt,\\
\end{array}
\end{equation}
where $y(\cdot)$ denote the uncertainties of demand
$$
y(s)=\displaystyle\frac{8}{3} s- W(s)
$$
and $u(\cdot)=\{u(s),0\leq s\leq 1\}$ is a control process taking values in a
compact set $U=[0,2]$.
Thus, we minimum  the following cost functional
\begin{equation}
J(u(\cdot))=E[X^u(\frac{1}{2})+X^u(1)], \label{1exa-1}%
\end{equation}
with the state constrains
$$
0\leq EX^u(\frac{1}{2}),\ \,0\leq  EX^u(1).
$$
Substituting  $X^u(\cdot)$ into equation (\ref{1exa-1}), one obtain
$$
J(u(\cdot))=E[2\displaystyle\int_0^{\frac{1}{2}}(u(t)-\frac{8}{3}t)dt+
\displaystyle\int_{\frac{1}{2}}^1(u(t)-\frac{8}{3}t)dt]
$$
and we can verify that
 \begin{equation}
   ( \bar{u}(t),\bar{X}(t))=
   \begin{cases}
    (\frac{8}{3}t,\displaystyle \int_0^tW(s)ds) &\mbox{ $0\leq t \leq \frac{1}{2}$,}\\
  (2,2t-\frac{4}{3}t^2-\frac{2}{3}+\displaystyle \int_0^tW(s)ds)&\mbox{ $ \frac{1}{2}  < t \leq 1$, }
   \end{cases}
  \end{equation}
is an optimal pair of systems (\ref{1exa-1}).

Recalling that it is difficult to get the adjoint equation for state process (\ref{3exp-1}). In order to get the related adjoint equations. We rewrite equation (\ref{3exp-1}) as follows,
$$
X^u(s)-W(s)s=\displaystyle \int_0^s(u(t)-\frac{8}{3}t)dt-\displaystyle\int_0^stdW(t).
$$
Denote that $\delta X^u(s)=X^u(s)-W(s)s$, thus
$$
d\delta X^u(s)=[u(s)-\frac{8}{3}s]ds-sdW(s)
$$
and we  have
 $$
 E[\delta X^u(\frac{1}{2})+\delta X^u(1)]=E[X^u(\frac{1}{2})+X^u(1)],
 $$
which means that $( \bar{u}(t),\bar{X}(t)-W(s)s)$ is the optimal pair of the following cost functional,
\begin{equation}
\label{3exc-2}
\delta J(u(\cdot))=E[\delta X^u(\frac{1}{2})+\delta X^u(1)],
\end{equation}
with the following constrains conditions,
$$
0\leq E\delta X^u(\frac{1}{2}),\ \,0\leq  E\delta X^u(1).
$$
Next, we introduce the following  first-order adjoint  equations for functional (\ref{3exc-2}).
\begin{equation*}%
\begin{array}
[c]{ll}%
d{p}(t)= q(t)dW(t),\quad \frac{1}{2}  < t < 1,\\
p(1)= -(\beta^0+\beta^2)
\end{array}
\label{prin-11}%
\end{equation*}
 and
 \begin{equation*}%
 \begin{array}
[c]{ll}%
d{p}(t)= q(t)dW(t),\quad  0\leq t < \frac{1}{2}, \\
p(\frac{1}{2})= -(\beta^0+\beta^1)+p(\frac{1}{2}^{+}).
\end{array}
\label{prin-12}%
\end{equation*}
where $(\beta^0,\beta^1,\beta^2)$ comes from Theorem \ref{ccc-th}. The solutions of first-order adjoint equations as follows:
 \begin{equation}
   ( p(t),q(t))=
   \begin{cases}
   (-(2\beta^0+\beta^1+\beta^2),0) &\mbox{ $ 0\leq t < \frac{1}{2}$, }\\
  (-(\beta^0+\beta^2),0)&\mbox{ $\frac{1}{2} < t < 1$.}
   \end{cases}
  \end{equation}
 Now, let $\beta^0+\beta^2\leq0$ and $2\beta^0+\beta^1+\beta^2= 0$, one obtain,
\begin{equation}
   H_u(\bar{X}(t)-W(t)t,\bar{u}(t),p(t),q(t))(v-\bar{u}(t))=
   \begin{cases}
    -(2\beta^0+\beta^1+\beta^2)(v-\frac{8}{3}t) &\mbox{ $ 0\leq t < \frac{1}{2}$, }\\
  -(\beta^0+\beta^2) (v-2) &\mbox{ $\frac{1}{2} < t \leq 1$,}\\

   \end{cases}
     \end{equation}
 with $v\in[0,2]$,  thus, the optimal control pair $(\bar{u}(\cdot),\bar{X}(\cdot))$ satisfies the Theorem \ref{ccc-th}.

 Notice that, under the same constraints for the  parameter $(\beta^0,\beta^1,\beta^2)$, we can verify other optimal pair for the model (\ref{3exc-2}).
\end{example}

\section{Necessary conditions for optimality}
In this section, we give the well known pontryagin's  stochastic maximum principle, which we will show the necessary conditions for optimal pairs.

Note that in cost functional (\ref{cost-1}),  we consider a multi-states  cost functional, which is very different from classical optimal control problem. Under a strong Fr\'{e}chet differentiable assumption, Yang \cite{Y16a,Y16b} studied the maximum principle for deterministic and stochastic systems. Also, we refer Gao and Yang \cite{GY16} for forward and backward stochastic system. In this study, we not only investigate the optimal control problem under a weak smooth condition, but also investigate a near optimal control model to cover \cite{Y16a,Y16b}. Since the optimal control set $U$ is convex, we only need to introduce the following first-order  adjoint equations:
\begin{equation}%
\begin{array}
[c]{ll}%
-d{p}(t)= & \{b_x(\bar{X}{(t)},\bar{u}(t))^{\mathrm{T}}p(t)+ \displaystyle\sum_{j=1}^d\sigma_x^j(\bar{X}{(t)},\bar{u}(t))^{\text{T}}q^j(t) \\
               &-f_x(\bar{X}{(t)},\bar{u}(t))\}dt-q(t)dW(t),\ t\in(t_{i},t_{i+1}),\\
p(t_{i+1})= &-E[\Psi_{x_{i+1}}(\bar{X}(t_1),\cdots,\bar{X}(t_n))|\mathcal{F}_{t_{i+1}}]+p(t_{i+1}^{+}),\text{ \ }i=0,1,\ldots,n-1,
\end{array}
\label{prin-1}%
\end{equation}
where "$\mathrm{T}$" means the transform of vector or matrix,  $t_{i+1}^{+}$\ is the right limit of $t_{i+1}$, $\Psi_{x_{i+1}}(x(t_1),\cdots,x(t_n))$ means the first derivative of $\Psi$ about $x(t_{i+1})$ and $p(t_{n}^{+})=0$, also, we set $t_0=0$.

Denote that
\begin{equation*}
H(x,u,p,q)=b(x,u)^{\text{T}}p+\sum_{j=1}^d\sigma^j (x,u)^{\text{T}}q^j-f(x,u),\text{ \  \ }%
(x,u,p,q)\in \mathbb{R}^m\times U\times \mathbb{R}^m\times \mathbb{R}^{m\times d}.%
\end{equation*}

The main result of this section is the following theorem.
\begin{theorem}
\label{Maximumprinciple} Let Assumptions (\ref{ass-b})-(\ref{ass-fai}) hold,
and $(\bar{u}(\cdot),\bar{X}(\cdot))$ be an optimal pair of (\ref{cost-2}).
Then there exists  $(p(\cdot),q(\cdot))$ satisfying the series of first-order adjoint equations (\ref{prin-1}) and respectively such that
\begin{equation}%
\begin{array}
[c]{ll}%
&H_u(\bar{X}(t),\bar{u}(t),p(t),q(t))(v-\bar{u}(t))\leq 0,\\
\end{array}
\label{prin-3}%
\end{equation}
for any $v\in U$ and $t \in(t_{i},t_{i+1})$, $i=0,1,\cdots,n-1$.
\end{theorem}

In the below, we will show the proof of the Theorem \ref{Maximumprinciple}.
The main difficult is to give the variational equation and adjoint equation
for multi-target terminal functional (\ref{cost-1}). For the limitation of technique in this section, we consider that $U$ is a convex set. Let $(\bar{u}(\cdot),\bar{X}(\cdot))$ be the given optimal pair. Let $0<\rho<1$  and $u(\cdot)+\bar{u}(\cdot)\in \mathcal{U}[0,T]$ be any given control. We define the following%
\[
u^{\rho}(t)=\bar{u}(t)+\rho u(t)=(1-\rho)\bar{u}(t)+\rho(u(t)+\bar{u}(t)),
\]
obviously, $u^{\rho}(\cdot)\in \mathcal{U}[0,T]$. The following Lemma is
useful for proving the Theorem \ref{Maximumprinciple}.

\begin{lemma}
\label{le-2} Let Assumptions (\ref{ass-b})-(\ref{ass-fai}) hold, and
$X^{\rho}(\cdot)$ be the solution of equation (\ref{ODE_1}) under the
control $u^{\rho}(\cdot)$, and $y^{}(\cdot)$  be the solutions
of the following equations:%
\begin{equation}
\label{apro-1}
\begin{array}
[c]{rl}%
d{y}(t)= & \big{[}b_{x}(\bar{X}{(t)},\bar{u}(t))y(t)+b_{u}(\bar{X}{(t)},\bar{u}(t)) u(t)\big{]}dt\\
&+ \displaystyle \sum_{j=1}^d\big{[} \sigma_{x}^j(\bar{X}{(t)},\bar{u}(t))y(t)+\sigma^j_u(\bar{X}{(t)},\bar{u}(t)) u(t)\big{]}dW^j(t), \\
y(0)= & 0,\quad t\in [0,T].
\end{array}
\end{equation}
Then%
\begin{equation}%
\begin{array}
[l]{l}%
\displaystyle\lim_{\rho\to 0}\displaystyle\sup_{t\in \lbrack0,T]}E\left|\rho^{-1} (X^{\rho}(t)-\bar{X}(t))-y(t)\right| =0,  \\
\end{array}
 \label{var-1}%
\end{equation}
and%
\begin{equation}%
\begin{array}
[c]{rl}
& \rho^{-1}J(u^{\rho}(\cdot))-J(\bar{u}(\cdot))\\
= &
{\displaystyle \sum \limits_{i=1}^{n}} E \big{[}\Psi_{x_i}(\bar{X}(t_1),\bar{X}(t_2),\cdots,\bar{X}(t_n))y(t_{i})\big{]} \\
& +
E{\displaystyle \int \limits_{0}^{T}}
\big{[}f_{x}(\bar{X}{(t)},\bar{u}(t))y(t)
+f_u(\bar{X}{(t)},\bar{u}(t)) u(t)\big{]}dt+o(1)\geq 0.
\end{array}
\label{var-3}%
\end{equation}

\end{lemma}
\textbf{Proof: }
Similar with the proof of Lemma 4.1 in \cite{P93}, we have equation (\ref{var-1}).
Note that%
\begin{equation}%
\begin{array}
[c]{rl}
& J(u^{\rho}(\cdot))-J(\bar{u}(\cdot))\\
= &E\big{[} \Psi(X^{\rho}(t_1),X^{\rho}(t_2),\cdots,X^{\rho}(t_n))-\Psi(\bar{X}(t_1),\bar{X}(t_2),\cdots,\bar{X}(t_n))\\
&+{\displaystyle \int \limits_{0}^{T}}
[f(X{}^{\rho}(t),u^{\rho}(t))-f(\bar{X}{(t)},\bar{u}(t))]dt\big{]},\\
\end{array}
\label{value-0}%
\end{equation}
which deduces that
\begin{equation}%
\begin{array}
[c]{rl}
& J(u^{\rho}(t))-J(\bar{u}(t))\\
= &\displaystyle \sum_{i=1}^nE\big{[} \Psi_{x_i}((\bar{X}(t_1),\bar{X}(t_2),\cdots,\bar{X}(t_n))({X}^{\rho}(t_i)-\bar{X}(t_i)) \big{]}\\
&+E\displaystyle \int \limits_{0}^{T}[f_x(\bar{X}{(t)},\bar{u}(t))(X^{\rho}(t)-\bar{X}{(t)})+
f_u(\bar{X}{(t)},\bar{u}(t))\rho u(t)]dt+o(\rho).
\end{array}
\label{value-1}%
\end{equation}
By equation (\ref{var-1}), it follows equation (\ref{var-3}).

This completes the proof. $\ \ \ \ \ \ \ \ \Box$
\bigskip

Based on the above Lemma, we now carry out the proof for Theorem
\ref{Maximumprinciple}.

\textbf{Proof of Theorem} \ref{Maximumprinciple}. For$\ t\in (
t_{i},t_{i+1}),$ applying the differential chain rule to $p(t)^{\rm{T}}y(t)$, we have%
\begin{equation}%
\begin{array}
[c]{rl}
&E\big{[} p(t_{i+1})^{\rm{T}}y(t_{i+1})-p(t_{i}^{+})^{\rm{T}}y(t_{i})\big{]}\\
=&E \big{[}-E[\Psi_{x_{i+1}}(\bar{X}(t_1),\cdots,\bar{X}(t_n))^{\rm T}|\mathcal{F}_{t_{i+1}}]y(t_{i+1})+p(t_{i+1}^{+})y(t_{i+1})-p(t_{i}^{+})y(t_{i}) \big{]}\\
=&E \big{[}-\Psi_{x_{i+1}}(\bar{X}(t_1),\cdots,\bar{X}(t_n))^{\rm T}y(t_{i+1})+p(t_{i+1}^{+})y(t_{i+1})-p(t_{i}^{+})y(t_{i}) \big{]}\\
=& E\displaystyle \int \limits_{t_i}^{t_{i+1}}\big{[}p(t)^{\text{T}}b_u(\bar{X}(t),\bar{u}(t))u(t)+\sum_{j=1}^d
q^j(t)^{\text{T}}\sigma_u^j (\bar{X}(t),\bar{u}(t))u(t)+f_x(\bar{X}(t),\bar{u}(t))^{\text{T}}y(t)   \big{]}dt

\end{array}
\label{max-1}%
\end{equation}
Adding by $i$ on the both sides of equation (\ref{max-1}) from $0$ to $n-1$, it follows
\begin{equation*}%
\begin{array}
[c]{rl}
&{\displaystyle \sum \limits_{i=1}^{n}} E \big{[}-\Psi_{x_i}(\bar{X}(t_1),\bar{X}(t_2),\cdots,\bar{X}(t_n))y(t_{i})-\displaystyle \int \limits_{t_{i-1}}^{t_{i}}
\big{[}f_{x}(\bar{X}{(t)},\bar{u}(t))^{\rm T}y(t)
+f_u(\bar{X}{(t)},\bar{u}(t))^{\rm T}u(t)\big{]}dt\\
=& {\displaystyle \sum \limits_{i=1}^{n}} E\displaystyle \int \limits_{t_{i-1}}^{t_{i}}\big{[}p(t)^{\text{T}}b_u(\bar{X}(t),\bar{u}(t))u(t)+\sum_{j=1}^d
q^j(t)^{\text{T}}\sigma_u^j (\bar{X}(t),\bar{u}(t))u(t)
-f_u(\bar{X}(t),\bar{u}(t))^{\text{T}}u(t)\big{]}dt
\end{array}
\end{equation*}
By Lemma \ref{le-2}, we obtain
\begin{equation*}%
\begin{array}
[c]{rl}
&{\displaystyle \sum \limits_{i=1}^{n}} E\displaystyle \int \limits_{t_{i-1}}^{t_{i}}\big{[}p(t)^{\text{T}}b_u(\bar{X}(t),\bar{u}(t))u(t)+\sum_{j=1}^d
q^j(t)^{\text{T}}\sigma_u^j (\bar{X}(t),\bar{u}(t))u(t)
-f_u(\bar{X}(t),\bar{u}(t))^{\text{T}}u(t)\big{]}dt\\
=&{\displaystyle \sum \limits_{i=1}^{n}} E\displaystyle \int \limits_{t_{i-1}}^{t_{i}}
\big{[} H_u(\bar{X}(t),\bar{u}(t),p(t),q(t)){u}(t)\big{]}dt\leq o(1),
\end{array}
\end{equation*}
Letting $\rho\to 0$, thus
\begin{equation*}%
\begin{array}
[c]{ll}%
&H_u(\bar{X}(t),\bar{u}(t),p(t),q(t))(v-\bar{u}(t))\leq 0,\\
\end{array}
\end{equation*}
for any $v\in U$ and $t \in(t_{i},t_{i+1})$, $i=0,1,\cdots,n-1$.

This completes the proof. $\ \ \ \ \ \ \ \ \Box$

\section{Sufficient conditions for optimality}
In this section, we consider another problem when $(\bar{u}(\cdot),\bar{X}(\cdot))$ is an optimal control pair of problem (\ref{cost-2}). Thus, we show the sufficient conditions for optimality in the following.
\begin{theorem}
\label{sufo}
Suppose Assumptions (\ref{ass-b})-(\ref{ass-fai}) hold and $\Psi(\cdot)$ is convex and $H(\cdot,\cdot,p(t),q(t))$ is concave for any $t\in(t_{i},t_{i+1})$ with $i=0,1,\cdots,n-1$ almost surely, and such that
\begin{equation}
\label{suf-1}
H(\bar{X}(t),\bar{u}(t),p(t),q(t))=\max_{u\in U}H(\bar{X}(t),u,p(t),q(t)),
\end{equation}
where $(p(\cdot),q(\cdot))$ is the solution of equation (\ref{prin-1}) with $(\bar{X}(\cdot),\bar{u}(\cdot))$. Thus, $(\bar{X}(\cdot),\bar{u}(\cdot))$ is an optimal pair of problem (\ref{cost-2}).
\end{theorem}

\noindent\textbf{Proof}: From the minimum condition (\ref{suf-1}), one obtain
$$
H_u(\bar{X}(t),\bar{u}(t),p(t),q(t))=0.
$$
Then for any given pair $(X^u(\cdot),u(\cdot))$ which solves equation (\ref{ODE_1}), and notice that $H(\cdot,\cdot,p(t),q(t))$ is concave, we have
\begin{equation}
\label{suf-2}
\begin{array}
[c]{ll}%
&\displaystyle\int_0^T\big{[} H(X^u(t),u(t),p(t),q(t))- H(\bar{X}(t),\bar{u}(t),p(t),q(t))\big{]}dt\\
\leq & \displaystyle\int_0^TH_x(\bar{X}(t),\bar{u}(t),p(t),q(t))(X^u(t)-\bar{X}(t)) dt.\\
\end{array}
\end{equation}
In the following, we introduce the approximation equation which $\delta X(t)=X^u(t)-\bar{X}(t) $,
\begin{equation}
\label{suf-3}
\begin{array}
[c]{rl}%
d\delta X(t)=&\big{[}b_x(\bar{X}(t),\bar{u}(t))\delta X(t)+\xi(t)  \big{]}dt\\
&+\displaystyle \sum_{j=1}^{d}\big{[} \sigma^j_x(\bar{X}(t),\bar{u}(t))\delta X(t)+\eta^j(t) \big{]}dW^j(t),\ t\in[0,T],\\
\delta X(t)=&0,
\end{array}
\end{equation}
where
\begin{equation*}
\begin{array}
[c]{rl}%
&\xi(t):=-b_x(\bar{X}(t),\bar{u}(t))\delta X(t)+b({X}^u(t),{u}(t))-b(\bar{X}(t),\bar{u}(t))\\
&\eta^j(t):=-\sigma^j_x(\bar{X}(t),\bar{u}(t))\delta X(t)+\sigma^j({X}^u(t),{u}(t))-\sigma^j(\bar{X}(t),\bar{u}(t)),\ 1\leq j \leq d. \\
\end{array}
\end{equation*}
For$\ t\in (
t_{i},t_{i+1}),$ applying the differential chain rule to $p(t)^{\rm{T}}\delta X(t)$, we have%
\begin{equation}%
\begin{array}
[c]{rl}
&E\big{[} p(t_{i+1})^{\rm{T}}\delta X(t_{i+1})-p(t_{i}^{+})^{\rm{T}}\delta X(t_{i})\big{]}\\
=&E \big{[}-E[\Psi_{x_{i+1}}(\bar{X}(t_1),\cdots,\bar{X}(t_n))^{\rm T}|\mathcal{F}_{t_{i+1}}]\delta X(t_{i+1})+p(t_{i+1}^{+})\delta X(t_{i+1})-p(t_{i}^{+})\delta X(t_{i}) \big{]}\\
=&E \big{[}-\Psi_{x_{i+1}}(\bar{X}(t_1),\cdots,\bar{X}(t_n))^{\rm T}\delta X(t_{i+1})+p(t_{i+1}^{+})\delta X(t_{i+1})-p(t_{i}^{+})\delta X(t_{i}) \big{]}\\
=& E\displaystyle \int \limits_{t_i}^{t_{i+1}}\big{[}p(t)^{\text{T}}\xi(t)+\sum_{j=1}^d
q^j(t)^{\text{T}}\eta^j(t)+f_x(\bar{X}(t),\bar{u}(t))^{\text{T}}\delta X(t)   \big{]}dt.
\end{array}
\label{max-1}%
\end{equation}
Adding by $i$ on the both sides of equation (\ref{max-1}) from $0$ to $n-1$, it follows
\begin{equation*}%
\begin{array}
[c]{rl}
&{\displaystyle \sum \limits_{i=1}^{n}} E \big{[}-\Psi_{x_i}(\bar{X}(t_1),\bar{X}(t_2),\cdots,\bar{X}(t_n))\delta X(t_{i})\big{]}\\
=& {\displaystyle \sum \limits_{i=1}^{n}} E\displaystyle \int \limits_{t_{i-1}}^{t_{i}}\big{[}p(t)^{\text{T}}\xi(t)+\sum_{j=1}^d
q^j(t)^{\text{T}}\eta^j(t)
+f_x(\bar{X}(t),\bar{u}(t))^{\text{T}}\delta X(t)\big{]}dt.
\end{array}
\end{equation*}
From the representation of $\xi$ and $\eta^j$, we have
\begin{equation*}%
\begin{array}
[c]{rl}
&{\displaystyle \sum \limits_{i=1}^{n}} E \big{[}-\Psi_{x_i}(\bar{X}(t_1),\bar{X}(t_2),\cdots,\bar{X}(t_n))\delta X(t_{i})\big{]}+E \displaystyle \int_0^T H_x(\bar{X}(t),\bar{u}(t),p(t),q(t))\delta X(t)dt\\
=&E\displaystyle \int_0^T\big{[}p(t)^{\rm T}(b({X}^u(t),{u}(t))-b(\bar{X}(t),\bar{u}(t)))+\displaystyle\sum\limits_{i=1}^{d}
q^j(t)^{\text{T}}(\sigma^j({X}^u(t),{u}(t))-\sigma^j(\bar{X}(t),\bar{u}(t)))\big{]}dt
\end{array}
\end{equation*}
By inequality (\ref{suf-2}), it follows that
\begin{equation}%
\label{suf-4}
\begin{array}
[c]{rl}
&{\displaystyle \sum \limits_{i=1}^{n}} E \big{[}-\Psi_{x_i}(\bar{X}(t_1),\bar{X}(t_2),\cdots,\bar{X}(t_n))\delta X(t_{i})\big{]}\\
\leq & E\displaystyle \int_0^T\big{[}f({X}^u(t),{u}(t))-f(\bar{X}(t),\bar{u}(t)) \big{]}dt.  \\
\end{array}
\end{equation}
Now, by the convexity of $\Psi$, we obtain
\begin{equation}%
\label{suf-5}
\begin{array}
[c]{rl}
&\displaystyle \sum \limits_{i=1}^{n}\Psi_{x_i}(\bar{X}(t_1),\bar{X}(t_2),\cdots,\bar{X}(t_n))\delta X(t_{i})\\
\leq & \Psi({X}^u(t_1),{X}^u(t_2),\cdots,{X}^u(t_n))
-\Psi(\bar{X}(t_1),\bar{X}(t_2),\cdots,\bar{X}(t_n)).
\end{array}
\end{equation}
Combining equations (\ref{suf-4}) and (\ref{suf-5}), we can verify that
$$
J(\bar{u}(\cdot))\leq J(u(\cdot)).
$$
By the arbitrary of $u(\cdot)$, we complete the proof. $\ \ \ \ \ \ \ \ \Box$

\section{Optimal control problem with state constraints}
In many applications of mathematics fiance and economics, we  may have different constraints at different time of the state process $X(\cdot)$, i.e., in the cost functional (\ref{cost-1}), let
$$
\Phi(X(t_1),X(t_2),\cdots,X(t_n))=\displaystyle\sum_{i=1}^n\phi(X(t_i)),
$$
then, in this case, our objective may be
$$
E\phi(X(t_i))\leq \alpha_i,\ \ i=1,2,\cdots,n,
$$
where $\alpha_i$ is a given constant. On the other hand, one may concern different combination of  $(X(t_1),X(t_2),\cdots,X(t_n))$, or a general case,
$$
E\Phi(X(t_1),X(t_2),\cdots,X(t_n))\leq \alpha,
$$
where $\alpha$ is a given constant.

In the following, we will first investigate the state equation (\ref{ODE_1}) with the below cost functional,
\begin{equation}
J(u(\cdot))=%
E\big{[}{\displaystyle \int \limits_{0}^{T}}
f(X{(t)},u(t))dt+\displaystyle\sum_{i=1}^n\phi(X(t_i))\big{]},\label{cost-3}%
\end{equation}
and the state process $X(\cdot)$ satisfies
\begin{equation}
\label{ccc-1}
E\phi(X(t_i))\leq \alpha_i,\ \ i=1,2,\cdots,n,
\end{equation}
where $\alpha_i$ is a given constant.

In order to prove the main result of this section, we introduce the following lemma which comes from Yong and Zhou Corollary 6.3 in \cite{Y99}.
\begin{lemma}
\label{ccc-le1}
Let $F:V\to \mathbb{R}$ be a continuous function on complete metric space $(V,\tilde{d})$. Given $\theta>0$ and $v_0\in V$ such that
$$
F(v_0)\leq \inf_{v\in V}F(v)+\theta.
$$
Then there exists a $v_{\theta}\in V$ such that
$$
F(v_{\theta})\leq F(v_0),\ \ \tilde{d}(v_{\theta},v_0)\leq \sqrt{\theta},
$$
and for all $v\in V$,
$$
- \sqrt{\theta}d(v_{\theta},v)\leq F(v)-F(v_{\theta}).
$$
\end{lemma}

Next, we present the main results of this section, the related Hamiltonian as follows,
\begin{equation*}
H(\beta^0,x,u,p,q)=b(x,u)^{\text{T}}p+\sum_{j=1}^d\sigma^j (x,u)^{\text{T}}q^j-\beta^0f(x,u),\text{ \  \ }
\end{equation*}
whith $(\beta^0,x,u,p,q)\in \mathbb{R}\times\mathbb{R}^m\times U\times \mathbb{R}^m\times \mathbb{R}^{m\times d}.$

\begin{theorem}
\label{ccc-th}
Let Assumptions (\ref{ass-b})-(\ref{ass-fai}) hold,
and $(\bar{u}(\cdot),\bar{X}(\cdot))$ be an optimal pair of (\ref{cost-3}).
Then there exists $(\beta^0,\beta^1,\cdots,\beta^n)\in \mathbb{R}^{n+1}$ satisfying
$$
\beta^0\geq 0,\ \ \left| \beta^0\right|^2+\displaystyle\sum_{j=1}^n\left| \beta^j\right|^2=1,
$$
and
$$
\displaystyle\sum_{j=1}^n\beta^j(\gamma^j-E\phi(\bar{X}(t_j)))\leq 0,\  \gamma^j\leq \alpha^j,\ j=1,2,\cdots,n,
$$
and the adapted solution $(p(\cdot),q(\cdot))$ satisfying the following series of first-order adjoint equations,
\begin{equation}%
\begin{array}
[c]{ll}%
-d{p}(t)= & \{b_x(\bar{X}{(t)},\bar{u}(t))^{\mathrm{T}}p(t)+ \displaystyle \sum_{j=1}^d\sigma_x^j(\bar{X}{(t)},\bar{u}(t))^{\text{T}}q^j(t) \\
               &-\beta^0f_x(\bar{X}{(t)},\bar{u}(t))\}dt-q(t)dW(t),\ t\in(t_{i-1},t_{i}),\\
p(t_{i})= &-(\beta^0+\beta^i)E[\Psi_{x_i}(\bar{X}(t_1),\cdots,\bar{X}(t_n))|\mathcal{F}_{t_i}]+p(t_{i}^{+}),\text{ \ }i=1,2,\ldots,n,
\end{array}
\label{prin-1}%
\end{equation}
and respectively such that
\begin{equation}%
\begin{array}
[c]{ll}%
&H_u(\beta^0,\bar{X}(t),\bar{u}(t),p(t),q(t))(v-\bar{u}(t))\leq 0,\\
\end{array}
\label{prin-3}%
\end{equation}
for any $v\in U$ and $t \in(t_{i},t_{i+1})$, $i=0,1,\cdots,n-1$.
\end{theorem}
\noindent\textbf{Proof}: Without loss of generality, we can assume that $J(\bar{u}(\cdot))=0$ where $(\bar{u}(\cdot),\bar{X}(\cdot))$ is the optimal pair of problem (\ref{cost-3}) with constraints (\ref{ccc-1}). For any $\theta>0$, we set
$$
J^{\theta}(u(\cdot))=\sqrt{\big{[}(J(u(\cdot))+\theta)^+\big{]}^2
+\displaystyle\sum_{i=1}^n\big{[}(E\phi(X^u(t_i))-\alpha_i)^+\big{]}^2}.
$$
From Assumption \ref{ass-fai}, one can verify that $J^{\theta}:\mathcal{U}[0,T]\to \mathbb{R}$ is continuous and satisfies
\begin{equation}
J^{\theta}(\bar{u}(\cdot))=\theta\leq \inf_{u\in\mathcal{U}[0,T]}J^{\theta}(u(\cdot))+\theta.
\end{equation}
Now, by Lemma \ref{ccc-le1}, there exists a $u^{\theta}(\cdot)\in \mathcal{U}[0,T]$ such that
\begin{equation}
\label{ccci-ev}
J^{\theta}(u^{\theta}(\cdot))\leq J^{\theta}(\bar{u}(\cdot))=\theta,\ \tilde{d}(u^{\theta}(\cdot),\bar{u}(\cdot))\leq \sqrt{\theta},
\end{equation}
where $\tilde{d}(u^1(\cdot),u^2(\cdot))=E\big{[}\displaystyle\int_0^T\left|u^1(t)-u^2(t)\right|^2dt\big{]}$. And we can check that $(\mathcal{U}[0,T],\tilde{d})$ is a complete metric space. Also, we have
\begin{equation*}
-\sqrt{\theta}\tilde{d}(u^{\theta}(\cdot),u(\cdot))\leq J^{\theta}(u(\cdot))-J^{\theta}(u^{\theta}(\cdot)),\ \forall u(\cdot) \in \mathcal{U}[0,T],
\end{equation*}
which deduces that
\begin{equation}
\label{cos-40}
J^{\theta}(u^{\theta}(\cdot))+\sqrt{\theta}\tilde{d}(u^{\theta}(\cdot),u^{\theta}(\cdot))\leq J^{\theta}(u(\cdot))+\sqrt{\theta}\tilde{d}(u^{\theta}(\cdot),u(\cdot)),\ \forall u(\cdot) \in \mathcal{U}[0,T].
\end{equation}
Thus, inequality (\ref{cos-40}) shows that $(u^{\theta}(\cdot),X^{\theta}(\cdot))$ is the optimal pair for the following cost functional
\begin{equation}
\label{cos-4}
J^{\theta}(u(\cdot))+\sqrt{\theta}\tilde{d}(u^{\theta}(\cdot),u(\cdot)),
\end{equation}
without the state constraint.

Since $U$ is a convex set, for any $\rho>0$, let $u^{\theta}(\cdot)+u(\cdot)\in \mathcal{U}[0,T]$, we define
$$
u^{\theta,\rho}(t)=u^{\theta}(t)+\rho u(t),
$$
which belongs to $ \mathcal{U}[0,T]$. It is easy to verify that
$$
\tilde{d}(u^{\theta,\rho}(\cdot),u^{\theta}(\cdot))=\rho E\displaystyle\int_0^T\left|u(t)\right|^2dt.
$$
For notation simplicity, we set $C^u=E\displaystyle\int_0^T\left|u(t)\right|^2dt$, by equation (\ref{cos-40}), one obtain
\begin{equation}%
\begin{array}
[c]{rl}%
-\sqrt{\theta}\rho C^u\leq & J^{\theta}(u^{\theta,\rho}(\cdot))-J^{\theta}(u^{\theta}(\cdot))\\
=&\displaystyle\frac{\big{[}(J(u^{\theta,\rho}(\cdot))+\theta)^+\big{]}^2-\big{[}(J(u^{\theta}(\cdot))+\theta)^+\big{]}^2
}
{J^{\theta}(u^{\theta,\rho}(\cdot))+J^{\theta}(u^{\theta}(\cdot))}\\
&+\displaystyle\frac{\sum_{j=1}^n\big{[}\big{[}(E\phi(X^{\theta,\rho}(t_j))-\alpha^j)^+\big{]}^2-
\big{[}(E\phi(X^{\theta}(t_j))-\alpha^j)^+\big{]}^2\big{]}}
{J^{\theta}(u^{\theta,\rho}(\cdot))+J^{\theta}(u^{\theta}(\cdot))},\\
\end{array}
\label{ccci-1}%
\end{equation}
where $X^{\theta,\rho}(\cdot))$ and $X^{\theta}(\cdot))$ are the related solution of equation (\ref{ODE_1}) with controls $u^{\theta,\rho}(\cdot)$ and $u^{\theta}(\cdot)$. Setting
\begin{equation}%
\begin{array}
[c]{ll}%
\beta^{0,\theta}=\displaystyle\frac{\big{[}J(u^{\theta}(\cdot))+\theta\big{]}^+}{J^{\theta}(u^{\theta}(\cdot))},\\
\beta^{j,\theta}=\displaystyle\frac{\big{[}E\phi(X^{\theta}(t_j))-\alpha^j\big{]}^+}
{J^{\theta}(u^{\theta}(\cdot))},\ j=1,2,\cdots,n.\\
\end{array}
\label{ccci-2}%
\end{equation}
Then, by the continuity of $J^{\theta}(\cdot)$ and Assumption \ref{ass-fai}, we have
\begin{equation}%
\begin{array}
[c]{rl}%
& J^{\theta}(u^{\theta,\rho}(\cdot))-J^{\theta}(u^{\theta}(\cdot))\\
=&\beta^{0,\theta}\big{[}J(u^{\theta,\rho}(\cdot))-J(u^{\theta}(\cdot))\big{]}+
\displaystyle\sum_{j=1}^n\beta^{j,\theta}\big{[}E\phi(X^{\theta,\rho}(t_j))
-E\phi(X^{\theta}(t_j))\big{]}+o(1),\\
=&E\big{[}\displaystyle\sum_{j=1}^n(\beta^{0,\theta}+\beta^{j,\theta})(\phi(X^{\theta,\rho}(t_j))
-\phi(X^{\theta}(t_j)))\\
&+\beta^{0,\theta}\displaystyle\int_0^T\big{[}f(X^{\theta,\rho}(t),u^{\theta,\rho}(t))-
f(X^{\theta}(t),u^{\theta}(t))\big{]}dt   \big{]}+o(\rho),
\end{array}
\label{ccci-3}%
\end{equation}
where $o(1)$ converges to $0$ when $\rho\to 0$.

Similar with Lemma \ref{le-2}, let $(\bar{X}(\cdot),\bar{u}(\cdot))$ be replaced by $(X^{\theta}(t),u^{\theta}(t))$, and $y(\cdot)$ be replaced by $\tilde{y}(\cdot)$ in equation (\ref{apro-1}). Thus, one obtain,
\begin{equation}%
\begin{array}
[c]{rl}
-\sqrt{\theta} C^u\leq &  \rho^{-1}\big{[}J^{\theta}(u^{\theta,\rho}(\cdot))-J^{\theta}(u^{\theta}(\cdot))\big{]}\\
\leq &
E \big{[}\displaystyle\sum_{j=1}^n(\beta^{0,\theta}+\beta^{j,\theta})\phi_{x}(X^{\theta}(t_j))
\tilde{y}(t_{j})\big{]} \\
& +
\beta^{j,\theta}E{\displaystyle \int \limits_{0}^{T}}
\big{[}f_{x}(X^{\theta}(t),u^{\theta}(t))\tilde{y}(t)
+f_u(X^{\theta}(t),u^{\theta}(t))u(t)\big{]}dt+o(1).
\end{array}
\label{ccci-4}%
\end{equation}
In addition, we introduce the following adjoint equation,
\begin{equation}%
\begin{array}
[c]{rl}%
-d{p}^{\theta}(t)= & \{b_x(X^{\theta}(t),u^{\theta}(t))^{\mathrm{T}}p^{\theta}(t)+ \sum_{j=1}^d\sigma_x^j(X^{\theta}(t),u^{\theta}(t))^{\text{T}}q^{j,{\theta}}(t) \\
               &-\beta^{0,\theta}f_x(X^{\theta}(t),u^{\theta}(t))\}dt-q^{\theta}(t)dW(t),\ t\in(t_{i-1},t_{i}),\\
p^{\theta}(t_{i})= &-(\beta^{0,\theta}+\beta^{i,\theta})E[\phi_{x}({X}^{\theta}(t_i))]+p(t_{i}^{+}),\text{ \ }i=1,\ldots,n,
\end{array}
\label{ccci-5}%
\end{equation}
where $q^{\theta}(\cdot)=(q^{1,\theta}(\cdot),q^{2,\theta}(\cdot),\cdots,q^{d,\theta}(\cdot))$.

Now, using the duality relation as in the proof of Theorem \ref{Maximumprinciple}, it follows that,
$$
{\displaystyle \sum \limits_{i=1}^{n}} E\displaystyle \int \limits_{t_{i-1}}^{t_{i}}
\big{[} H_u(\beta^{0,\theta},X^{\theta}(t),u^{\theta}(t),p^{\theta}(t),q^{\theta}(t)){u}(t)\big{]}dt\leq o(1)+\sqrt{\theta} C^u,
$$
Notice that $o(1) \to 0$ when $\rho\to 0$. Thus, letting $\rho\to 0$, one obtain
\begin{equation}
\label{ccci-6}
{\displaystyle \sum \limits_{i=1}^{n}} E\displaystyle \int \limits_{t_{i-1}}^{t_{i}}
\big{[} H_u(\beta^{0,\theta},X^{\theta}(t),u^{\theta}(t),p^{\theta}(t),q^{\theta}(t)){u}(t)\big{]}dt\leq \sqrt{\theta} C^u.
\end{equation}

From inequality (\ref{ccci-ev}), it follows that $u^{\theta}(\cdot)$ converges  to $\bar{u}(\cdot)$ under $\tilde{d}$ as $\theta\to 0$. Then, by Assumptions \ref{ass-b}, \ref{assb-b2} and \ref{ass-fai}, and basic theory of stochastic differential equation, we have
$$
\displaystyle \sup_{0\leq t\leq T}E\left|X^{\theta}(t)-\bar{X}(t)\right| \to 0,
$$
as $\theta\to 0$. By equation (\ref{ccci-2}), we have
\begin{equation}
\label{unq-1}
\left|\beta^{0,\theta}\right|^2+\displaystyle\sum_{j=1}^n\left|\beta^{j,\theta}\right|^2=1.
\end{equation}
Thus, we can choice a sequence $\{\theta_k\}_{k=1}^{\infty}$ satisfying $\displaystyle\lim_{k\to\infty}\theta_k=0$ and such that the limitations of $\beta^{0,\theta_k}$ and $\beta^{j,\theta_k}$ exist and we set
\begin{equation}
\begin{array}
[c]{ll}
\beta^{0}=\displaystyle\lim_{k\to\infty}\beta^{0,\theta_k},\\
\beta^{j}=\displaystyle\lim_{k\to\infty}\beta^{j,\theta_k},\\
\end{array}
\end{equation}
with $j=1,2,\cdots,n$. From equation (\ref{unq-1}), we have
$$
\left|\beta^{0}\right|^2+\displaystyle\sum_{j=1}^n\left|\beta^{j}\right|^2=1,
$$
and
$$
\displaystyle\sum_{j=1}^n\beta^j(\gamma^j-E\phi(\bar{X}(t_i)))\leq 0,\  \gamma^j\leq \alpha^j,\ j=1,2,\cdots,n.
$$
Similarly, we can prove that
$$
\displaystyle \sup_{0\leq t\leq T}E\big{[}\left|p^{\theta_k}(t)-p(t)\right|^2+\int_0^T\left|q^{\theta_k}(t)-q(t)\right|^2 \big{]} dt\to 0,
$$
as $k\to\infty$.  Letting $k\to \infty$, from equation (\ref{ccci-6}), we have
\begin{equation}
{\displaystyle \sum \limits_{i=1}^{n}} E\displaystyle \int \limits_{t_{i-1}}^{t_{i}}
\big{[} H_u(\beta^{0},\bar{X}(t),\bar{u}(t),p(t),q(t)){u}(t)\big{]}dt\leq 0.
\end{equation}
Thus, we complete this proof. $\ \ \ \ \ \ \ \ \ \ \ \  \Box$
\begin{remark}
\label{re-2}
Similarly with the proof in Theorem \ref{ccc-th}, we can deal with other constraints conditions, i.e.,
$$
\underline{\alpha}_i\leq E\phi(X(t_i))\leq \overline{\alpha}_i,\ \ i=1,2,\cdots,n,
$$
where $(\underline{\alpha}_i,\overline{\alpha}_i)_{i=1}^n$ are  given constants, or
$$
\underline{\alpha}\leq E\Phi(X(t_1),X(t_2),\cdots,X(t_n))\leq \overline{\alpha}.
$$

\end{remark}

\section{Near optimality for general case}
Recalling that in our previous paper \cite{Y16b}, we consider the state process (\ref{ODE_1}) with the following general cost functional,
\begin{equation}
J(u(\cdot))=%
E\big{[}{\displaystyle \int \limits_{0}^{T}}
f(X{(t)},u(t))dt+\Phi(X_{[0,T]})\big{]},\label{gcosf-1}%
\end{equation}
where $X_{[0,T]}:=X(s)_{0\leq s \leq T},$ which is the path of $X(\cdot)$ from $0$ to $T$. In \cite{Y16b}, under a strong assumption about  Fr\'{e}chet derivatives, the maximum principle for cost functional (\ref{gcosf-1}) is given by solving a sequence of new adjoint equations. In this section, we will remove the strong assumption about Fr\'{e}chet derivatives and develop a near maximum principle for the cost functional (\ref{gcosf-1}) via the argument in Section 3. For notation simplicity, we set $m=d=1$.

\begin{assumption}
\label{assn-1}
Suppose $\Psi$ is Lipschatiz continuous on $\mathbb{C}[0,T]$, there exists a constant $c>0$ such that
$$
\left|\Phi(x^1_{[0,T]})-\Phi(x^2_{[0,T]})\right|\leq c \max_{0\leq t\leq T}\left|x^1(t)-x^2(t)\right|,
$$
where $x^1_{[0,T]},x^2_{[0,T]}\in \mathbb{C}[0,T]$, and $\mathbb{C}[0,T]$ is the set of continuous functions over $[0,T]$.
\end{assumption}

By Assumptions \ref{assn-1}, one obtain that there exists a larger enough integer $N>0$, for $n>N$ such that
$$
\left|\Phi(x(t_1),x(t_2),\cdots,x(t_n))-\Phi(x_{[0,T]})\right|\leq c\max_{1\leq j\leq n}\sup_{t_{j-1}\leq t\leq t_{j}} \left|x(t)-x(t_j)\right|,
$$
with $t_0=0$ and $c$ is the constant in Assumptions \ref{assn-1}. Next, we define the approximation function for $\Phi(x(t_1),x(t_2),\cdots,x(t_n))$ as follows,
\begin{equation}
\label{gae-1}
\begin{array}
[c]{rl}
{\Phi}^{\varepsilon}(x(t_1),x(t_2),\cdots,x(t_n))=&\displaystyle\int_{\mathbb{R}^n}\big{[}\Phi(x(t_1),x(t_2),\cdots,x(t_n))
\psi^{\varepsilon}(y_1-x(t_1))\\
&\times\psi^{\varepsilon}(y_2-x(t_2))\cdots
\psi^{\varepsilon}(y_n-x(t_n))\big{]}dy_1dy_2\cdots dy_n,\\
\end{array}
\end{equation}
with $(x(t_1),x(t_2),\cdots,x(t_n))\in \mathbb{R}^{n}$, and $\psi^{\varepsilon}(x)=\frac{1}{\sqrt{2\pi \varepsilon^2}}e^{-\frac{x^2}{2\varepsilon^2}}$ for $x\in \mathbb{R}$. Thus, we have the following Lemma.
\begin{lemma}
\label{visle-1}
There exists a constant $C>0$ such that
\begin{equation*}
\begin{array}
[c]{ll}
\left|{\Phi}^{\varepsilon}(x(t_1),x(t_2),\cdots,x(t_n))-{\Phi}(x(t_1),x(t_2),\cdots,x(t_n))\right|\leq C\varepsilon,\\
\end{array}
\end{equation*}
$\forall (x(t_1),x(t_2),\cdots,x(t_n))\in \mathbb{R}^{n}$.
\end{lemma}
\textbf{Proof}: We just prove the case $n=1$. For general case, we can use the same method. Similarly, we can obtain the other inequalities. By equation (\ref{gae-1}), for fixed $x$, we have,
\begin{equation}
\begin{array}
[c]{ll}
&|{\Phi}^{\varepsilon}(x)-\Phi(x)|\\
\leq & \displaystyle\int_{\mathbb{R}}\left|\Phi(x^0)-\Phi(x)
\right|\frac{1}{\sqrt{2\pi \varepsilon^2}}e^{-\frac{(x^0-x)^2}{2\varepsilon^2}}dx^0\\
\leq & C \displaystyle\int_{\mathbb{R}}\left|x^0-x
\right|\frac{1}{\sqrt{2\pi \varepsilon^2}}e^{-\frac{(x^0-x)^2}{2\varepsilon^2}}dx^0\\
=&C\varepsilon\displaystyle\int_{\mathbb{R}}{|\tilde{x}^0|}\frac{1}{\sqrt{2\pi}}
e^{-\frac{(\tilde{x}^0)^2}{2}}d\tilde{x}^0\\
\leq & C \varepsilon,
\end{array}
\end{equation}
where $C$ will change line by line.

This completes the proof. $\ \ \ \ \ \ \ \ \  \Box$

\begin{remark}
\label{re-1}
Notice that $\psi\in\mathbb{C}^{\infty}[\mathbb{R}]$, by the property of convolution, we obtain that $\Phi^{\varepsilon}(x(t_1),x(t_2),\cdots,x(t_n))$ is second differentiable about $(x(t_1),x(t_2),\cdots,x(t_n))\in\mathbb{R}^n$.
\end{remark}

In the following, we introduce the near optimal control problem,
\begin{equation}
{J}^{\varepsilon}(u(\cdot))=%
E\big{[}{\displaystyle \int \limits_{0}^{T}}
f(X{(t)},u(t))dt+\Phi^{\varepsilon}(X(t_1),X(t_2),\cdots,X(t_n))\big{]}.\label{gcosf-2}%
\end{equation}
By Assumptions \ref{ass-b} and \ref{assb-b2}, we can obtain the following results.
\begin{theorem}
\label{gle-1}
Let Assumptions \ref{ass-b} and \ref{assb-b2} hold, then, there exists a constant $C>0$ such that
$$
\left|\displaystyle\inf_{u\in\mathcal{U}[0,T]}J(u(\cdot))-\displaystyle\inf_{u\in\mathcal{U}[0,T]}{J}^{\varepsilon}(u(\cdot))\right|\leq C\varepsilon
$$
\end{theorem}

Now, we can use the results in the above sections to investigate the near optimal control theory for optimal control problem (\ref{gcosf-1}).

\end{document}